\documentclass[10pt]{amsart}
\usepackage{amsmath}
\usepackage[usenames,dvipsnames]{color}
\usepackage{parskip}
\usepackage{amsfonts}
\usepackage{amscd}
\usepackage[centertags]{amsmath}
\usepackage{amssymb}
\usepackage[all,cmtip]{xy}
\usepackage[english]{babel}
\usepackage{tabularx}
\usepackage{amsxtra}
\usepackage{euscript}
\usepackage[T1]{fontenc}
\usepackage{doc, exscale, fontenc, latexsym, syntonly}
\usepackage{amsfonts}
\usepackage{amsthm}
\usepackage{graphicx}

\newtheorem{theorem}{Theorem}[section]
\newtheorem{corollary}[theorem]{Corollary}
\newtheorem{lemma}[theorem]{Lemma}
\newtheorem{proposition}[theorem]{Proposition}
\theoremstyle{definition}
\newtheorem{definition}[theorem]{Definition}
\newtheorem{example}[theorem]{Example}

\numberwithin{figure}{section}
\numberwithin{table}{section}

\newcommand*\acknowledgment[1]{%
	\begingroup\noindent
	\rightskip\leftskip
	\begin{flushleft}\textbf{\large Acknowledgment.}\, #1%
		\par\vspace*{1mm}\end{flushleft}\endgroup}
\begin{document}

\title[COUNTEREXAMPLES FOR TOPOLOGICAL COMPLEXITY IN DIGITAL IMAGES]{COUNTEREXAMPLES FOR TOPOLOGICAL COMPLEXITY IN DIGITAL IMAGES}

\author{MEL\.{I}H \.{I}S and \.{I}SMET KARACA}
\date{\today}

\address{\textsc{Melih Is}
Ege University\\
Faculty of Sciences\\
Department of Mathematics\\
Izmir, Turkey}
\email{melih.is@ege.edu.tr}
\address{\textsc{Ismet Karaca}
Ege University\\
Faculty of Science\\
Department of Mathematics\\
Izmir, Turkey}
\email{ismet.karaca@ege.edu.tr}



\begin{abstract}
Digital topology has its own working conditions and sometimes differs from the normal topology. In the area of topological robotics, we have important counterexamples in this study to emphasize this red line between a digital image and a topological space. We indicate that the results on topological complexities of certain path-connected topological spaces show alterations in digital images. We also give a result about the digital topological complexity number using the genus of a digital surface in discrete geometry.
\end{abstract}

\maketitle

\section{Introduction}
\quad Digital topology makes it possible to analyze digital images by transferring topological properties into itself. This analyzing process provides advantages in technologies that include especially computer science and image analysis. After that digital topology is introduced by Rosenfeld \cite{Ros:1979} with its most primitive concepts, it has been involved in the studies of a wide variety of fields. One of the most significant of these studies points to the subject of robotics, including algebraic topology over the past 20 years.

\quad Farber \cite{Farber:2003} assigns a positive integer called the topological complexity number $TC(X)$ for each path-connected topological space $X$ and makes inferences about the complexity of the motion area in which a robot is located with obstacles. Afterwards, as the structure of the required topological space differs, the results about computing the topological complexity number also get vary. For example, if $X$ is a connected Lie group, then the topological complexity number $TC(X)$ of $X$ equals the Lusternik Schnirelmann category $cat(X)$ of $X$ \cite{Farber:2004}. Davis \cite{Davis:2017}, \cite{Davis2:2017} has different approaches on the topological complexity number of certain path-connected topological spaces such as a circle or a Klein bottle. The higher topological complexity $TC_{n}$, defined by Rudyak \cite{Rudyak:2010}, is the Schwarz genus of a special fibration and a positive integer such as $cat$ and $TC$. These are all homotopy invariants. Many essential concepts of the algebraic topology, for instance cohomological cup product \cite{Farber:2003}, help to calculate these numbers. Hence, it is the main task to have concrete ideas about the topological complexities of the digital images by making more use of the algebraic topology methods. See \cite{Farber:2008} for more information about the concepts, results and methods for topological robotics. 

\quad Digital topology has started to incorporate the algebraic topology methods in time and this process has not ended yet. One of the richest content of these methods is cohomology rings of topological spaces. Ege and Karaca \cite{EgeKaraca:2013} presented the general framework related to the digital simplicial cohomology groups and introduced the cohomological cup-products in digital images. In digital topology, researchers generally deal with the two-dimensional or three-dimensional digital images and they call them as digital curves and digital surfaces, respectively. Rosenfeld \cite{Ros:1970} defined the digital curve and stated that the digital image $X$ is said to be a simple (closed) curve if each point of $X$ has exactly two adjacent points in $X$. On the other hand, the definition of the digital surface was given by Morgenthaler and Rosenfeld \cite{MorgenRosen:1981}. Digital surfaces are so important not only for digital topology but also for discrete and computational geometry. Chen \cite{Chen:2004} used one of the most important notions related to topology and geometry in digital images for the first time: digital manifolds. A digital manifold can be regarded as a discretization of a manifold. In other saying, it is a combinatorial manifold which is defined in digital images. Chen and Rong \cite{ChenRong:2010} improved a powerful method for computing genus and the Betti numbers of a digital image.

\quad In digital topology, not everything is the same as in ordinary topology. As an example, K{\"u}nneth Theorem does not work for digital images \cite{EgeKaraca:2013}. This leads to us to show that cohomological cup product does not hold for digital images \cite{MelihKaraca}. Therefore, the digital interpretations of $cat$ \cite{BorVer:2018}, $TC$ \cite{KaracaIs:2018} and $TC_{n}$ \cite{MelihKaraca} have different results from the topological spaces. In this study, we focus on such results and give counterexamples in digital topology. First, we remind the basic notions of the digital topology, give some concepts and results from discrete geometry and mention about certain definitions of the main tools of robotics in Section \ref{section2}. Thus, we recall a few deliberately selected results on mostly Farber's topological complexity numbers. The main goal of this study gives counterexamples for each properties that is valid in the usual topology but not in digital topology. This shows that digital topology cannot be thought as the same as the usual topology. They often have different impact areas in mathematics because the methods of algebraic topology can lose its influence such as in the example of K{\"u}nneth Theorem and this is the strength of digital topology conditions. 

\quad In the first example, we consider the digital image $X = [0,1]_{\mathbb{Z}}$. Since the closed interval $[0,1]$ is contractible in topological spaces, the higher topological complexity of $[0,1]$ for any $n$ is $1$. Moreover, Proposition \ref{p1} supports this result using the cohomological cup-product but in digital images, using the digital cup-product does not give the same result. Hence, the digital image $X = [0,1]_{\mathbb{Z}}$ is a counterexample of that Proposition \ref{p1} does not hold for digital images. Second, we show that the diagonal map of digital images does not coincide with the digital cup product homomorphism. In the next example, we choose a special digital image that the Betti number is $2$. In the usual topology, it is expected that the topological complexity number of a space is $3$, when $b_{1}(X) = 2$ but we prove that it is possible that the digital topological complexity number of such an image is less than $3$. We present a result about the digitally connected digital curves with considering the first Betti numbers. Our next counterexample is constructed on showing that the digital version of the topological complexity number of the wedge of $2$ sphere $S^{2}$ does not have to be $3$. Finally, we study with two digital surfaces such that genus of them are $1$ and $2$. In topology, the topological complexity number of a compact orientable surface of genus $1$ is $3$ whereas the topological complexity number of a compact orientable surface of genus $2$ is $5$. We prove that both of this results does not hold for digital surfaces. We show that the digital topological complexity number of digital simple closed surfaces of genus $0$, $1$ and $2$ are $1$, $2$ and $3$, respectively. 

\section{Preliminaries}\label{section2}
\quad This section basically consists of the fundamental notions of the digital topology. In addition, remarkable points of the study of topological robotics are mentioned in the section.  

\quad Let $\mathbb{Z}^{m}$ be a set of only integers points of Euclidean space $\mathbb{R}^{m}$. Let $X$ be any finite subset of $\mathbb{Z}^{m}$. Then $(X,\kappa)$ is a \textit{digital image} \cite{Boxer:1999}, where $\kappa$ is an adjacency relation for the points of $X$. Let $x$ and $y$ be any two different points in $\mathbb{Z}^{m}$. Let $k$ be a positive integer such that $k$ is less than or equal to $m$. $x$ and $y$ are \textit{$c_{k}-$adjacent} \cite{Boxer:1999} if there are at most $k$ indices $i$ such that $|x_{i} - y_{i}| = 1$ and for all other indices $i$ such that $|x_{i} - y_{i}| \neq 1$, $x_{i} = y_{i}$. Whereas two points of $\mathbb{Z}$ has only $c_{1} = 2-$adjacency, two points of $\mathbb{Z}^{2}$ has $c_{1} = 4-$adjacency and $c_{2} = 8-$adjacency. Similarly, two points of $\mathbb{Z}^{3}$ has $c_{1} = 6-$adjacency, $c_{2} = 18-$adjacency and $c_{3} = 26-$adjacency. 

\quad Let $X$ be a digital image in $\mathbb{Z}^{m}$. $X$ is \textit{$\kappa-$connected} \cite{Herman:1993} if and only if for every pair of different points $x,y \in X$, there is a set $\{x_{0},x_{1}, ...,x_{l}\}$ of points of $X$ such that $x=x_{0}$, $y=x_{l}$ and $x_{i}$ and $x_{i+1}$ are $\kappa-$adjacent where $i=0,1,...,l-1$. Let $f : X_{1} \rightarrow X_{2}$ be a digital map such that $X_{1}$ has an adjacency relation $\kappa_{1}$ in $\mathbb{Z}^{m_{1}}$ and $X_{2}$ has an adjacency relation $\kappa_{2}$ in $\mathbb{Z}^{m_{2}}$. Then $f$ is \textit{$(\kappa_{1},\kappa_{2})-$continuous} \cite{Boxer:1999} if, for any $\kappa_{1}-$connected subset $U_{1}$ of $X_{1}$, $f(U_{1})$ is $\kappa_{2}-$connected. Moreover, if $f$ is both digitally continuous and bijective and the digital map $f^{-1}$ is $(\kappa_{2},\kappa_{1})-$continuous, then $f$ is \textit{$(\kappa_{1},\kappa_{2})-$isomorphism} \cite{Boxer2:2006}.

\quad Let $X$ and $Y$ be any two digital images such that $f,f^{'} : X \rightarrow Y$ are $(\kappa_{1},\kappa_{2})-$continuous maps. If there exists a positive integer $n$ and a digital map \[G : X \times [0,n]_{\mathbb{Z}} \rightarrow Y\] such that the following conditions hold, then $f$ and $f^{'}$ are said to be \textit{digitally $(\kappa_{1},\kappa_{2})-$homotopic} \cite{Boxer:1999} in $Y$. $G$ is called a \textit{digital $(\kappa_{1},\kappa_{2})-$homotopy} between $f$ and $f^{'}$. 

\begin{itemize}
	\item for all $x \in X$, $G(x,0)=f(x)$ and $G(x,n) = f^{'}(x)$;
	
	\item for all $x \in X$, the digital map $G_{x} : [0,n]_{\mathbb{Z}} \rightarrow Y$, defined by $G_{x}(t) = G(x,t)$, is $(2,\kappa_{2})-$continuous, for all $t \in [0,n]_{\mathbb{Z}}$;
	
	\item for all $t \in [0,n]_{\mathbb{Z}}$, the digital map $G_{t} : X \rightarrow Y$, defined by $G_{t}(x) = G(x,t)$, is $(\kappa_{1},\kappa_{2})-$continuous, for all $x \in X$.
\end{itemize} 

\quad Let $X$ be a digital image. If the identity map $X \rightarrow X$ is $(\kappa,\kappa)-$homotopic to a constant map \[c : X \longrightarrow X\]
\hspace*{5.9cm}$x \longmapsto c(x) = c_{0}$, \\
for all $x \in X$, then $(X,\kappa)$ is \textit{$\kappa-$contractible} \cite{Boxer:1999}.

\quad Let $X_{1}$ be a digital image with $\lambda_{1}-$adjacency and let $X_{2}$ be a digital image with $\lambda_{2}-$adjacency. Given two points $(x_{1},x_{2})$ and $(x_{1}^{'},x_{2}^{'})$ in the cartesian product digital image $X_{1} \times X_{2}$. Then \textit{$(x_{1},x_{2})$ and $(x_{1}^{'},x_{2}^{'})$ are adjacent in $X_{1} \times X_{2}$} \cite{BoxKar:2012} if one of the following conditions holds:
\begin{itemize}
	\item $x_{1} = x_{1}^{'}$ and $x_{2} = x_{2}^{'}$; or
	
	\item $x_{1} = x_{1}^{'}$ and $x_{2}$ and $x_{2}^{'}$ are $\lambda_{2}-$adjacent; or
	
	\item $x_{1}$ and $x_{1}^{'}$ are $\lambda_{1}-$adjacent and $x_{2} = x_{2}^{'}$; or	
	
	\item $x_{1}$ and $x_{1}^{'}$ are $\lambda_{1}-$adjacent and $x_{2}$ and $x_{2}^{'}$ are $\lambda_{2}-$adjacent.
\end{itemize}

\quad Let $X$ be a digital image. If $f : [0,n]_{\mathbb{Z}} \rightarrow X$ is a $(2,\kappa)-$continuous map such that $f(0)=x$ and $f(n)=x^{'}$, then $f$ is called a \textit{digital path} \cite{Boxer:2006} from $x$ to $x^{'}$ in $X$. A simple closed $\kappa-$curve of $r \geq 4$ points in a digital image $X$ is a sequence $\{g(0), g(1), ... , g(r-1)\}$ of images of the $\kappa-$path $g : [0,r-1]_{\mathbb{Z}} \rightarrow X$ such that $g(m)$ and $g(n)$ are \textit{$\kappa-$adjacent} if and only if $n = (m \pm 1)mod \hspace*{0.2cm}r$ \cite{Boxer:2005}. A \textit{digital surface} \cite{MorgenRosen:1981} is the set of surface points each of which has two adjacent components not in the surface in its neighborhood. Let $\kappa$ be an adjacency relation defined on $\mathbb{Z}^{m}$. A \textit{$\kappa-$neighbor} \cite{Herman:1993} of $x \in \mathbb{Z}^{m}$ is a point in $\mathbb{Z}^{m}$ that is $\kappa-$adjacent to $x$.

\begin{theorem}\cite{ChenRong:2010}
	If $X$ is a closed digital surface, then the genus of the surface is  
	\[g = 1 + \frac{(|M_{5}|+2.|M_{6}|-|M_{3}|)}{8},\]
	
	where $M_{i}$ is a set of points with $i-$neighbors.
\end{theorem}

\quad Let $(X,\kappa)$ be a digital image in $\mathbb{Z}^{2}$. \textit{The digital wedge union $X \vee X$} \cite{Han:2005} is the disjoint union of two $X$ with only one point $x_{0}$ in common and for any different elements $x$ and $y$ in $X$, $x$ and $y$ are not $\kappa-$adjacent to each other except the point $x_{0}$.

\quad Let $PX$ denote the set of all digitally continuous paths $\alpha : [0,n]_{\mathbb{Z}} \rightarrow X$ in $X$. $\pi : PX \rightarrow X \times X$ is a digitally continuous map that takes any digitally continuous path $\alpha$ in $X$ to the pair of its starting and ending points $(\alpha(0),\alpha(n))$.

\begin{definition}\cite{KaracaIs:2018}
	\textit{Digital topological complexity number $TC(X,\kappa)$} is the minimal number $l$ such that $U_{1}, U_{2}, ..., U_{l}$ is a cover of $X \times X$ and for all $1 \leq i \leq l$, there is a digitally continuous map $s_{i} : U_{i} \rightarrow PX$ such that $\pi \circ s_{i} = id_{U_{i}}$. If no such $l$ exists we will set $TC(X,\kappa) = \infty$.
\end{definition}

\quad In the definition of the digital topological complexity number, the digital continuity of $s_{i}$ is required. Due to this reason, the task is to define an adjacency relation between two digital paths. Let $\alpha_{1} : [0,n_{1}]_{\mathbb{Z}} \rightarrow X$ and $\alpha_{2} : [0,n_{2}]_{\mathbb{Z}} \rightarrow X$ be any digitally continuous paths in $X$. Then $\alpha_{1}$ and $\alpha_{2}$ are $\lambda-$connected on $PX$, if for all $t$ times, $\alpha_{1}(t)$ and $\alpha_{2}(t)$ are $\lambda-$connected. Note that $\alpha_{1}$ and $\alpha_{2}$ can have different $t$ times. Without loss of generality, assume that $n_{1} < n_{2}$. In this case, the steps of the shortest path are extended such that $\alpha_{1}(n_{1} + l) = \alpha_{1}(n_{1})$, where $0 \leq l \leq n_{2}-n_{1}$. Hence, $\alpha_{1}$ is synchronized with $\alpha_{2}$. 

\begin{definition}\cite{MelihKaraca}
	Let $p : X \longrightarrow X^{'}$ be a digital fibration. \textit{The digital Schwarz genus of $p$} is defined as the minimum number $l$ such that $U_{1}, U_{2}, ..., U_{l}$ is a cover of $X^{'}$ and for all $1 \leq j \leq l$, there is a digitally continuous map $t_{j} : U_{j} \rightarrow X$, such that $p \circ t_{j} = id_{U_{j}}$.
\end{definition} 

\begin{definition}\cite{MelihKaraca}
	Let $X$ be any $\kappa$-connected digital image. Let $J_{n}$ be the wedge of $n-$digital intervals $[0,m_{1}]_{\mathbb{Z}}, [0,m_{2}]_{\mathbb{Z}}, ... , [0,m_{n}]_{\mathbb{Z}}$, for a positive integer $n$, where $0_{i} \in [0,m_{i}]$, $i = 1, ... ,n$, are identified. Then \textit{the digital higher topological complexity $TC_{n}(X,\kappa)$} is defined by the digital Schwarz genus of the digital fibration $$e_{n}:X^{J_{n}} \rightarrow X^{n}$$
	\hspace*{5.7cm} $f \longmapsto (f(m_{1})_{1},...,f(m_{n})_{n})$,
	
	where $(m_{i})_{k}$, $k=1, ..., n$ denotes the endpoints of the $i-$th interval for each $i$. 
\end{definition}

\quad Note that for $n=2$, the digital higher topological complexity number coincides with the topological complexity number \cite{MelihKaraca}. The digital higher topological complexity is a homotopy invariant of digital images. $TC_{n}$ is a natural lower bound for $TC_{n+1}$.

\begin{definition}\cite{BorVer:2018}
	Let $X$ be a digital image. Then \textit{the digital Lusternik-Schnirelmann category $cat_{\kappa}(X)$} of $X$ is defined to be the minimum number $l$ such that there is a cover $U_{1}, U_{2}, ..., U_{l}$ of $X$ such that for all $i = 1, ..., l$, each $U_{i}$ is $\kappa-$contractible to a point in $X$.
\end{definition}

\begin{definition}\cite{MelihKaraca2}
	Let $X$ be a digital image with $\kappa-$adjacency and let $(X,\circ)$ be a group. Let $X \times X$ be equipped the minimum of adjacency relation for the cartesian product. $(X,\kappa,\ast)$ is a \textit{$\kappa-$topological group} if the digital maps
	\[\alpha : X \times X \longrightarrow X \hspace*{1.0cm} \text{and} \hspace*{1.0cm} \beta : X \to X,\] 
	\hspace*{3.5cm}$(x,x^{'}) \longmapsto x \circ x^{'}$ \hspace*{2.4cm} $x \longmapsto x^{-1}$, 
	
	are digitally continuous, for all $x$, $x^{'} \in X$.
\end{definition}

\begin{theorem}\label{t1}\cite{MelihKaraca2}
	Let $(H,\kappa,\circ)$ be a $\kappa-$topological group such that $(H,\kappa)$ is digitally connected and $n > 1$. Then \[TC_{n}(H,\kappa) = cat_{\kappa_{\ast}}(H^{n-1}),\]
	where $\kappa_{\ast}$ is an adjacency relation for $H^{n-1}$.
\end{theorem}

\quad We use the cohomology groups of the digital image $MSS^{'}_{6}$ \cite{Han:2007} in the next chapter, so we completely need to know these values. As a special example of the computation of cohomology groups of digital images, Proposition \ref{p2} answers it. Let us recall the example with its proof. (See \cite{Gulseli:2014} and \cite{GulseliKaraca:2017} for more examples of digital curves and digital surfaces about computations in details of both homology and cohomology groups of the digital images.)

\begin{proposition}\cite{Gulseli:2014}\label{p2}
	Let $MSS^{'}_{6}$ be a digital image which consists of $8$ points $p_{0}$, $p_{1}$, $p_{2}$, $p_{3}$, $p_{4}$, $p_{5}$, $p_{6}$ and $p_{7}$ in $\mathbb{Z}^{3}$, where 
	\begin{eqnarray*}
		&&p_{0} = (1,0,0), p_{1} = (1,1,0), p_{2} = (1,1,1), p_{3} = (1,0,1), \\ 
		&&p_{4} = (0,0,1), p_{5} = (0,1,1), p_{6} = (0,1,0), p_{7} = (0,0,0).
	\end{eqnarray*}
    The digital cohomology groups of $MSS^{'}_{6}$ are \[H^{q,6}(MSS^{'}_{6}) = \begin{cases} \mathbb{Z}, & q=0 \\ \mathbb{Z}^{5}, & q=1 \\ 0 & q \neq 0,1. \end{cases}\]
\end{proposition}

\begin{proof}
	Let $p_{7} < p_{4} < p_{6} < p_{5} < p_{0} < p_{3} < p_{1} < p_{2}$. $C_{0}^{6}(MSS^{'}_{6})$ and $C_{1}^{6}(MSS^{'}_{6})$ are free abelian groups with bases $0-$simplexes $<p_{0}>, ..., <p_{7}>$ and $1-$simplexes
	\begin{eqnarray*}
	&&e_{0} = <p_{0}p_{1}>, e_{1} = <p_{0}p_{3}>, e_{2} = <p_{1}p_{2}>, e_{3} = <p_{6}p_{1}>\\
	&&e_{4} = <p_{5}p_{2}>, e_{5} = <p_{3}p_{2}>, e_{6} = <p_{4}p_{3}>, e_{7} = <p_{7}p_{4}>\\
	&&e_{8} = <p_{4}p_{5}>, e_{9} = <p_{6}p_{5}>, e_{10} = <p_{7}p_{6}>, e_{11} = <p_{7},p_{0}>,
	\end{eqnarray*}
    respectively. For $q > 1$,  $C_{q}^{6}(MSS^{'}_{6}) = \{0\}$. Then the short exact sequence
    \begin{eqnarray*}
    	0 \stackrel{\partial_{2}}{\longrightarrow}C_{1}^{6}(MSS^{'}_{6}) \stackrel{\partial_{1}}{\longrightarrow}C_{0}^{6}(MSS^{'}_{6}) \stackrel{\partial_{0}}{\longrightarrow} 0
    \end{eqnarray*}
    is obtained. Since 
    \begin{eqnarray*}
    	&&C^{0,6}(MSS^{'}_{6}) = Hom(C_{0}^{6}(MSS^{'}_{6}),\mathbb{Z}) \hspace*{0.3cm} \text{and}\\
    	&&C^{1,6}(MSS^{'}_{6}) = Hom(C_{1}^{6}(MSS^{'}_{6}),\mathbb{Z}),
    \end{eqnarray*} 
    we have the sequence
    \begin{eqnarray*}
    	0 \stackrel{\delta^{-1}}{\longrightarrow}C^{0,6}(MSS^{'}_{6}) \stackrel{\delta^{0}}{\longrightarrow}C^{1,6}(MSS^{'}_{6}) \stackrel{\delta^{1}}{\longrightarrow} 0.
    \end{eqnarray*}
    It is easy to see that
    \begin{eqnarray*}
    	&&\partial_{1}(e_{0}) = p_{1}-p_{0}, \hspace*{0.2cm} \partial_{1}(e_{1}) = p_{3}-p_{0},\hspace*{0.2cm}  \partial_{1}(e_{2}) = p_{2}-p_{1},\hspace*{0.2cm}  \partial_{1}(e_{3}) = p_{1}-p_{6}, \\
    	&&\partial_{1}(e_{4}) = p_{2}-p_{5},\hspace*{0.2cm}  \partial_{1}(e_{5}) = p_{2}-p_{3},\hspace*{0.2cm}  \partial_{1}(e_{6}) = p_{3}-p_{4},\hspace*{0.2cm}  \partial_{1}(e_{7}) = p_{4}-p_{7}, \\
    	&&\partial_{1}(e_{8}) = p_{5}-p_{4},\hspace*{0.2cm}  \partial_{1}(e_{9}) = p_{5}-p_{6},\hspace*{0.2cm}  \partial_{1}(e_{10}) = p_{6}-p_{7},\hspace*{0.2cm}  \partial_{1}(e_{11}) = p_{0}-p_{7}. 
    \end{eqnarray*}
    So we find $0-$cochains
    \begin{eqnarray*}
    	&&\delta^{0}p_{0}^{*} = -e_{0} - e_{1} + e_{11}, \hspace*{0.2cm} \delta^{0}p_{1}^{*} = e_{0} - e_{2} + e_{3}, \hspace*{0.2cm} \delta^{0}p_{2}^{*} = e_{2} + e_{4} + e_{5}, \\
    	&&\delta^{0}p_{3}^{*} = e_{1} - e_{5} + e_{6}, \hspace*{0.2cm} \delta^{0}p_{4}^{*} = -e_{6} + e_{7} - e_{8}, \hspace*{0.2cm} \delta^{0}p_{5}^{*} = -e_{4} + e_{8} + e_{9}, \\
    	&&\delta^{0}p_{6}^{*} = -e_{3} - e_{9} + e_{10}, \hspace*{0.2cm} \delta^{0}p_{7}^{*} = -e_{7} - e_{10} - e_{11}.
    \end{eqnarray*}
    Moreover, we get
    \begin{eqnarray*}
    	\delta^{0}(\displaystyle \sum_{i=1}^{7}n_{i}p_{i}^{*}) &=& e_{0}(-n_{0} + n_{1}) + e_{1}(-n_{0} + n_{3}) + e_{2}(-n_{1} + n_{2}) + e_{3}(n_{1} - n_{6})\\
    	&+& e_{4}(n_{2} - n_{5}) + e_{5}(n_{2} - n_{3}) + e_{6}(n_{3} - n_{4}) + e_{7}(n_{4} - n_{7})\\
    	&+& e_{8}(-n_{4} + n_{5}) + e_{9}(n_{5} - n_{6}) + e_{10}(n_{6} - n_{7}) + e_{11}(n_{0} - n_{7}).
    \end{eqnarray*}
    If $\delta^{0}(\displaystyle \sum_{i=1}^{7}n_{i}p_{i}^{*}) = 0$, then we find $n_{0} = n_{1} = ... = n_{7} = n$. This shows that $Ker\delta^{0} = \mathbb{Z}$. Therefore, we have that $Z^{0,6}(MSS^{'}_{6}) \cong \mathbb{Z}$. Since $Im\delta^{-1} = \{0\}$, we have that $B^{0,6}(MSS^{'}_{6}) \cong \{0\}$. Hence, we get $H^{0,6}(MSS^{'}_{6}) \cong \mathbb{Z}$. In addition, we have that 
    \begin{eqnarray*}
    B^{1,6}(MSS^{'}_{6}) = Im\delta^{0} &=& \{t_{0}e_{0} + t_{1}e_{1} + t_{2}e_{2} + t_{3}e_{3} + t_{4}e_{4} \\
    &+& (t_{2} - t_{1})e_{5}, + t_{5}e_{6} + t_{6}e_{7} + (-t_{1} + t_{2} - t_{4} + t_{5})e_{8} \\
    &+& (t_{2} + t_{3} - t_{4})e_{9} + (t_{0} - t_{1} - t_{3} + t_{5} + t_{6})e_{10} \\
    &+& (- t_{1} + t_{5} + t_{6})e_{11} \hspace*{0.2cm}| \hspace*{0.2cm} i = 0,1,...,6, \hspace*{0.2cm} \forall t_{i} \in \mathbb{Z}\} \cong \mathbb{Z}^{7}.
    \end{eqnarray*}
    By the definition, we see that $Z^{1,6}(MSS^{'}_{6}) = Ker\delta^{1} \cong \mathbb{Z}^{12}$. This implies that $H^{1,6}(MSS^{'}_{6}) = Z^{1,6}(MSS^{'}_{6}) / B^{1,6}(MSS^{'}_{6}) \cong \mathbb{Z}^{5}$. Finally, we get our result
    \begin{eqnarray*}
    	H^{q,6}(MSS^{'}_{6}) = \begin{cases} \mathbb{Z}, & q=0 \\ \mathbb{Z}^{5}, & q=1 \\ 0 & q \neq 0,1. \end{cases}
    \end{eqnarray*}
\end{proof}

\quad Before ending the section, we list some results on computing the topological complexity numbers or the higher topological complexity numbers in topological spaces. Each of the results is invalid in digital images. The counterexamples of the digital images are exhibited in Section \ref{section3}.

\begin{proposition}\cite{Rudyak:2010}\label{p1}
	Let $\Delta_{n} : X \rightarrow X^{n}$ be the diagonal map. If there exist $v_{i} \in H^{\ast}(X^{n};M_{i})$, $i = 1, ... ,k$, for which $(\Delta_{n})_{\ast}v_{i} = 0$ and $v_{1} \smile v_{2} \smile ... \smile v_{k} \neq 0$ in the cohomology $H^{\ast}(X^{n};M_{1} \otimes M_{2} \otimes ... \otimes M_{k})$, then $TC_{n}(X) \geq k+1$. 
\end{proposition}

\begin{theorem}\cite{Farber:2004} \label{t2}
	Let $X$ be a connected graph. Then 
	\[TC(X) = 
	\begin{cases}
	1, & \text{if} \hspace*{0.2cm}b_{1}(X) = 0, \\ 2, & \text{if} \hspace*{0.2cm} b_{1}(X) = 1, \\ 3, & \text{if} \hspace*{0.2cm} b_{1}(X) \geq 2,
	\end{cases}\]
	where $b_{1}(X)$ denotes the first Betti number of $X$.
\end{theorem}

\begin{lemma}\label{l1}\cite{Farber:2004}
	Let $X$ denote the wedge of $k$ spheres $S^{n}$. Then
	\[TC(X) = 
	\begin{cases}
	2, & \text{if} \hspace*{0.2cm} k=1 \hspace*{0.2cm} \text{and} \hspace*{0.2cm} n \hspace*{0.1cm} \text{is} \hspace*{0.1cm} \text{odd}, \\ 3, & \text{if} \hspace*{0.2cm} \text{either} \hspace*{0.2cm} k>1, \hspace*{0.2cm} \text{or} \hspace*{0.2cm} n \hspace*{0.1cm} \text{is} \hspace*{0.1cm} \text{even}.
	\end{cases}\]
\end{lemma}

\begin{example}\label{e2}\cite{Farber:2006}
	Let $X = \sum_{g}$ be a compact orientable surface of genus $g$. Then
	\[ TC(X) = \begin{cases}
	3, & \text{if} \hspace*{0.2cm} g \leq 1 \\
	5, & \text{if} \hspace*{0.2cm} g > 1.
	\end{cases}\] 
\end{example}
\section{MAIN RESULTS}\label{section3}
\quad The aim of this section is to emphasize the differences between the topological spaces and the digital images in the sense of the objects of topological robotics. The results and the examples in this section wish to indicate that the studies of the topological complexity numbers are so valuable and can be transfered into the digital images to get interesting results. This does not lose the value of the studies, contrarily it enriches the studies by carrying into the different platform. The results can be instructive for the future study of the topological complexity in digital image processing.

\quad First, we show that Proposition \ref{p1} is not valid in digital images. 

\begin{example}\label{e1}
	Take $n=3$ and $X = [0,1]_{\mathbb{Z}}$ in the Proposition \ref{p1}. Let $M$ be a field. The diagonal map $\Delta_{3} : (X,2) \rightarrow (MSS^{'}_{6},6)$ induces the digital homomorphism on digital cohomology with the coefficent $\mathbb{M}$ for the first dimension \[(\Delta_{3})_{\ast} : H^{1,6}(MSS^{'}_{6};\mathbb{M}) \rightarrow H^{1,2}(X;\mathbb{M}).\]
	Proposition \ref{p2} tells us that $(\Delta_{3})_{\ast}$ is basically the digital map $\mathbb{M}^{5} \rightarrow 0$. This shows that $ker((\Delta_{3})_{\ast}) = \mathbb{M}^{5}$. For any nonzero element $(v_{1}, v_{2}, v_{3}, v_{4}, v_{5})$ of $\mathbb{M}^{5}$, we have that $(\Delta_{3})_{\ast} (v_{1}, v_{2}, v_{3}, v_{4}, v_{5}) = 0$. Hence, we obtain $k \geq 1$. On the other hand, we shall show that $TC_{3}(X,2) = 1$. $([0,1]_{\mathbb{Z}},2,\circ)$ is a $2-$topological group, where $\circ$ is defined by
	\[\circ (a,b) = \begin{cases}
	0, & a=b \\ 1, & a \neq b,
	\end{cases}\]
	for all $a$, $b \in [0,1]_{\mathbb{Z}}$. By Theorem \ref{t1}, $TC_{3}(X,2) = cat_{4}(X^{2})$. Since $X^{2}$ is $4-$contractible, we obtain $TC_{3}(X,2) = 1$. As a result, $TC_{3}(X,2) < k+1$.
\end{example}

 \quad As a result of this example, [Proposition 3.2,\cite{GrantLuptonOprea:2013}] becomes meaningless in digital images because the diagonal map $\Delta_{n}$ cannot be identified with the digital cup product homomorphism $\cup_{n}(X)$ in [Definition 3.1,\cite{GrantLuptonOprea:2013}]. We shall explain it with an example in digital images:

\begin{example}
	Let $\mathbb{M}$ be a field. Consider the digital cohomology induced by the diagonal map $\Delta_{3}$, with the coefficient $M$ for the first dimension,  given in Example \ref{e1}. Although $(\Delta_{3})_{\ast}$ is the digital map $\mathbb{M}^{5} \rightarrow 0$, the digital cup product homomorphism \[\smile_{3} : H^{1,2}(X;\mathbb{M}) \otimes_{\mathbb{M}} H^{1,2}(X;\mathbb{M}) \otimes_{\mathbb{M}} H^{1,2}(X;\mathbb{M}) \rightarrow H^{3,2}(X;\mathbb{M})\]
	is the digital map $0 \otimes_{\mathbb{M}} 0 \otimes_{\mathbb{M}} 0 \rightarrow 0$. This shows that $Ker(\smile_{3}) = 0$. Hence, $Ker(\smile_{3})$ does not coincidence with $Ker((\Delta_{3})_{\ast})$. Thus, we conclude that \[nil(Ker(\smile_{3})) < TC_{3}(X,2) < nil(Ker(\Delta_{3})_{\ast}).\]
\end{example}

\quad In the next example we show that Theorem \ref{t2} does not hold for digital images for the case $b_{1}(X) \geq 2$.

\begin{figure*}[h]
	\centering
	\includegraphics[width=0.40\textwidth]{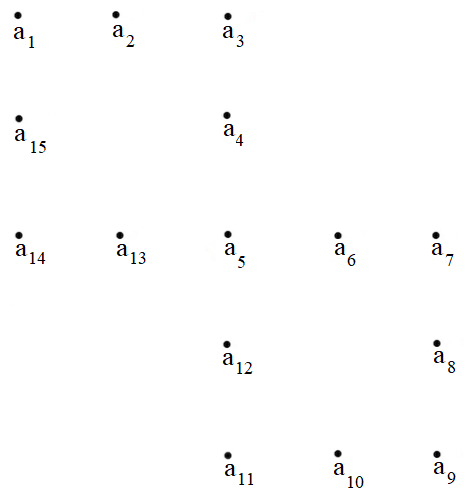}
	\caption{The digital image $X$.}
	\label{fig1:figure1}
\end{figure*}
\begin{example}\label{e3}
	Let $X$ be a digital image as shown in the Figure \ref{fig1:figure1}. The first Betti number of $X$ is $2$, because $X$ has two digital quadrilateral holes. On the other hand, we shall show that $TC(X,4) = 2$. $X$ is not $4-$contractible, so $TC(X,4) > 1$. Let $\alpha = \{a_{1}, a_{15}, a_{14}, a_{13}, a_{5}, a_{6}, a_{7}, a_{8}, a_{9}\}$ and $\beta = \{a_{1}, a_{2}, a_{3}, a_{4}, a_{5}, a_{12}, a_{11}, a_{10}, a_{9}\}$ be the subsets of $X$ as shown in the Figure \ref{fig2:figure2}.
	\begin{figure*}[h]
		\centering
		\includegraphics[width=0.90\textwidth]{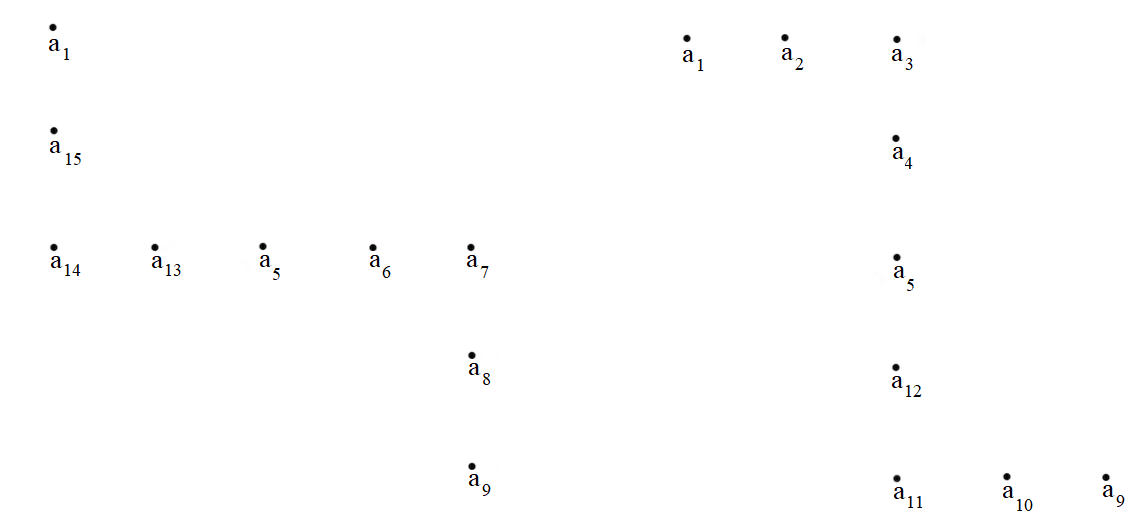}
		\caption{The digital images $\alpha$ on the left and $\beta$ on the right in $X$.}
		\label{fig2:figure2}
	\end{figure*}
     $X \times X$ can be written as the union of the following sets
    \begin{eqnarray*}
    	&&U_{1} = \{(x,y) \in X \times X \hspace*{0.2cm}|\hspace*{0.2cm} (x,y) \in \beta \hspace*{0.2cm} \text{or} \hspace*{0.2cm} (x,y) \in \beta \times \alpha \hspace*{0.2cm} \text{or} \hspace*{0.2cm} (x,y) \in \alpha \times \beta\} \\
    	&&\hspace*{4.0cm}\text{and}\\
    	&&U_{2} = \{(x,y) \in X \times X \hspace*{0.2cm}|\hspace*{0.2cm} (x,y) \in \alpha\}.
    \end{eqnarray*}
    By the definition, the minimal number is $2$ for the digital topological complexity number of $X$. Finally, we conclude that $TC(X,4) = 2$.
\end{example}

\begin{theorem}
	Let $X$ be a digitally connected digital curve and $b_{1}(X)$ denote the first Betti number of $X$. Then 
	\[TC(X) = 
	\begin{cases}
	1, & \text{if} \hspace*{0.2cm}b_{1}(X) = 0, \\ 2, & \text{if} \hspace*{0.2cm} b_{1}(X) \geq 1.
	\end{cases}\]
\end{theorem}
\begin{figure*}[h]
	\centering
	\includegraphics[width=1.00\textwidth]{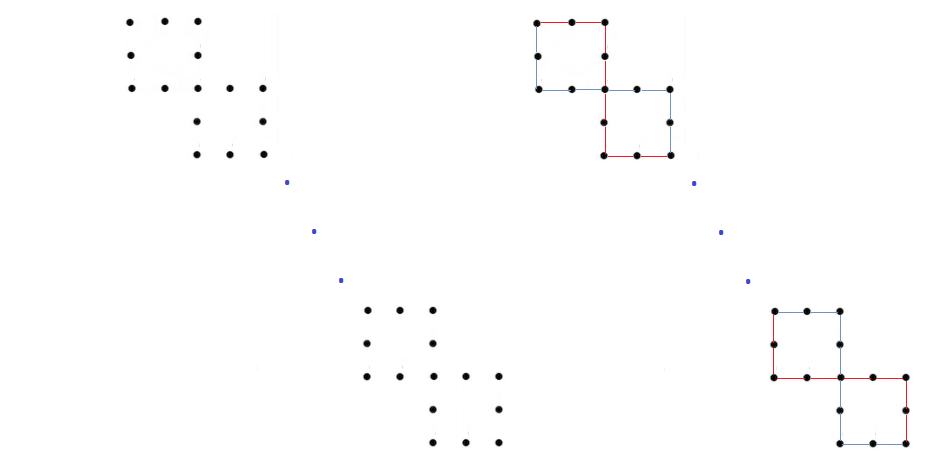}
	\caption{The left one is the digital image with the Betti number $n > 1$. The red and blue lines describe the paths passing through the elements of $T_{1}$ and $T_{2}$, respectively.}
	\label{fig10:figure10}
\end{figure*}
\begin{proof}
	Let $X$ be a digitally connected digital curve. If $b_{1}(X) = 0$, then $X$ has no quadrilateral holes in it. This means that $X$ is $k-$contractible. Hence, we find that $TC(X,\kappa) = 1.$ If $b_{1}(X) = 1$, there is only one hole in $X$ and we cover $X$ by two sets such like that in Example \ref{e3}. The digital homotopy invariance property of the digital topological complexity number concludes that $TC(X,\kappa) = 2$. Consider the case $b_{1}(X) > 1.$ Let $b_{1}(X) = n$, where $n > 1$. Then $X$ is covered by two sets $T_{1}$ and $T_{2}$, where $T_{i} = W_{i} \times W_{i}$ (see Figure \ref{fig10:figure10} for the images $W_{i}$), for $i=1,2$. By the homotopy invariance property, we get $TC(X,\kappa) = 2$, for $b_{1}(X)>1$.
\end{proof}

\quad The following example shows that Lemma \ref{l1} is not true in the digital meaning of topological complexity.

\begin{example}
	 Take both the numbers $k$ and $n$ as $2$ in Lemma \ref{l1}. Then we have the space $S^{2} \vee S^{2}$. The digitally equivalent of this space is digitally isomorphic to the digital image $X = MSC^{'}_{6} \vee MSC^{'}_{6}$ with the wedge point $(0,0,0)$ as shown in Figure \ref{fig3:figure3}. Since the digital topological complexity is a digital homotopy invariant, we observe that $TC(MSC^{'}_{6} \vee MSC^{'}_{6},6) = 3$. On the other hand, $MSC^{'}_{6} \vee MSC^{'}_{6}$ is $6-$contractible. Indeed, the homotopy between the identity map on $MSC^{'}_{6} \vee MSC^{'}_{6}$ and the constant map at the point $(0,0,0)$ has $3$ steps. For $t \in [0,2]_{\mathbb{Z}}$, we explain how the contraction homotopy
	 \begin{eqnarray*}
	 	H : (MSC^{'}_{6} \vee MSC^{'}_{6}) \times [0,3]_{\mathbb{Z}} \rightarrow MSC^{'}_{6} \vee MSC^{'}_{6}
	 \end{eqnarray*}
     works (see Figure \ref{fig4:figure4}). For the time $t=0$, each point of $MSC^{'}_{6} \vee MSC^{'}_{6}$ is located in its original place. For the time $t=1$, the upper $MSC^{'}_{6}$ is deformed to the square whose points are $b_{1}$, $b_{2}$, $b_{3}$ and $(0,0,0)$. Similarly, for the lower $MSC^{'}_{6}$, we have the square whose points are $b_{4}$, $b_{5}$, $b_{6}$ and $(0,0,0)$. For the time $t=2$, the squares are deformed to the point $(0,0,0)$. Consequently, this is a contradiction because the $6-$contractibility of the image gives us that $TC(MSC^{'}_{6} \vee MSC^{'}_{6},6) = 1$.
	 \begin{figure*}[h]
	 	\centering
	 	\includegraphics[width=0.40\textwidth]{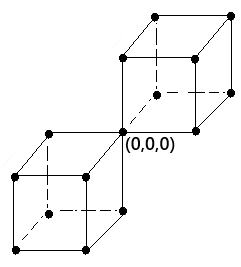}
	 	\caption{The digital image $MSC^{'}_{6} \vee MSC^{'}_{6}$ and the wedge point $(0,0,0)$.}
	 	\label{fig3:figure3}
	 \end{figure*}  
     \begin{figure*}[h]
     	\centering
     	\includegraphics[width=0.70\textwidth]{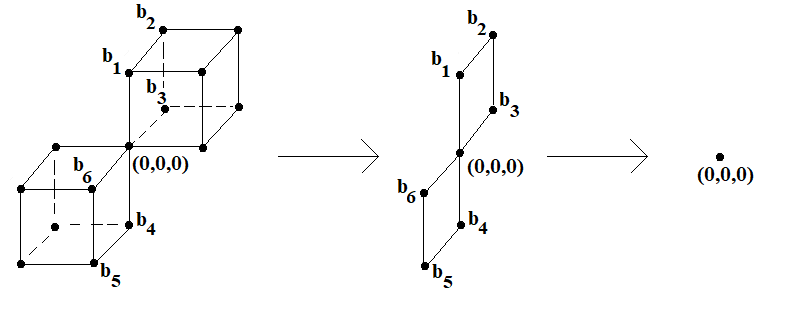}
     	\caption{The digital contraction of $MSC^{'}_{6} \vee MSC^{'}_{6}$.}
     	\label{fig4:figure4}
     \end{figure*}	
\end{example}

\begin{corollary}
	Let $(X,\kappa)$ be a digital image that consists of the wedge of $r$ digital cube in $\mathbb{Z}^{n}$, for $n>1$. Then
	\[TC(X,\kappa) = 
	\begin{cases}
	1, & \text{if} \hspace*{0.2cm} n=3 \hspace*{0.2cm} \text{or} \hspace*{0.2cm} n=2 \hspace*{0.2cm} \text{and} \hspace*{0.2cm} \kappa = 8 \\ 2, & \text{if} \hspace*{0.2cm} n=2 \hspace*{0.2cm} \text{and} \hspace*{0.2cm} \kappa = 4 \\  
	\end{cases}\]
\end{corollary}

\begin{proof}
	If $n=3$ and $r=1$, then the topological complexity number of $X$ is the same as the topological complexity number of the digital image $MSC^{'}_{6}$. Since this is a $6-$contractible space, the topological complexity of $X$ is independence from the number $r$ and the adjacencies in $\mathbb{Z}^{3}$. If $n=2$, then the topological complexity number of $X$ is equal to the topological complexity number of $MSC_{4}$. By $TC(MSC_{4},4) = 2$ and $TC(MSC_{4},8) = 1$, the result holds.
\end{proof}

\quad The next result shows that Example \ref{e2} is not true for digital images.

\begin{example}\label{e4}
	Let $X$ and $Y$ be two digital images as in Figure \ref{fig7:figure7}.
	\begin{figure*}[h]
		\centering
		\includegraphics[width=0.70\textwidth]{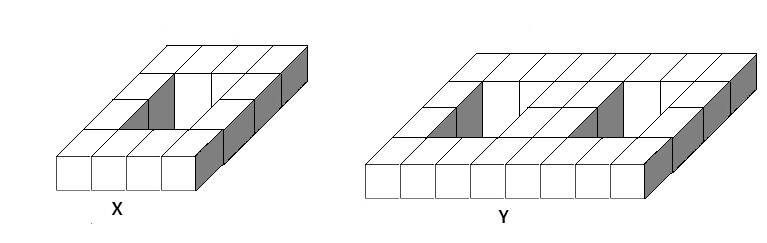}
		\caption{The digital images with genus $1$ and genus $2$, respectively.}
		\label{fig7:figure7}
	\end{figure*}
    Chen \cite{Chen:2004} shows that $X$ is an example of the digital closed surface with genus $1$ and $Y$ is an example of the digital closed surface with genus $2$. By Example \ref{e2}, the digital topological complexity number of $X$ is $3$ and the topological complexity number of $Y$ is $5$. However, we shall show this is a contradiction.
    \begin{figure*}[h]
    	\centering
    	\includegraphics[width=0.50\textwidth]{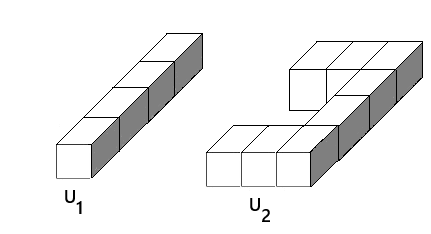}
    	\caption{Divide the digital image $X$ into $2$ parts.}
    	\label{fig8:figure8}
    \end{figure*}
    Take $U_{1}$ and $U_{2}$ as in Figure \ref{fig8:figure8}. If we define $V_{i} = U_{i} \times U_{i}$, for all $i=1,2$, then we get $TC(X,6) = 2$.
    Similarly, take $T_{1}$, $T_{2}$ and $T_{3}$ as in Figure \ref{fig9:figure9}. If we define $W_{i} = T_{i} \times T_{i}$, for all $1 \leq i \leq 3$, then we obtain that $TC(Y,6) = 3$.
    \begin{figure*}[h]
    	\centering
    	\includegraphics[width=0.70\textwidth]{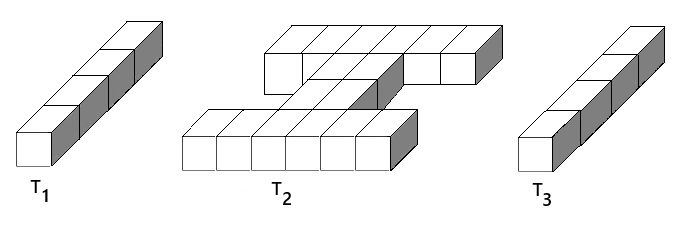}
    	\caption{Divide the digital image $Y$ into $3$ parts.}
    	\label{fig9:figure9}
    \end{figure*}
\end{example}

\begin{corollary}
	Let $X$ be a digital simple closed surface of genus $g$ with $6-$adjacency. Then 
	\[TC(X,6) = \begin{cases}
	1, & \text{if} \hspace*{0.2cm} g= 0 \\
	2, & \text{if} \hspace*{0.2cm} g = 1 \\
	3, & \text{if} \hspace*{0.2cm} g = 2.
	\end{cases}\]
\end{corollary}

\begin{proof}
	It is enough to show that $TC(X,6) = 1$, when $g=0$. The image $MSC^{'}_{6}$ is an example of a digital surface of genus $0$. It is clear that $TC(MSC^{'}_{6},6) = 1$. The fact that the digital topological complexity number is a homotopy invariance of digital images gives the desired result.
\end{proof}

\section{Conclusion}
\quad The subjects of robotics are widely studied in every aspect of science. Algebraic topology has put these studies on a different and effective ground by means of configuration spaces. Our purpose in the future is to use digital topology instead of ordinary topology and obtain remarkable results in the development of computer sciences. The evaluation of robotics' works in digital topology is still very new. First, the differences and the similarities should be determined. Therefore, determining the differences between the topological complexity number of a topological space and the digital topological complexity number of a digital image is the main goal of this study. 

\quad We begin with considering how the relation between cohomological cup product homomorphism and the diagonal map works for digital images. Second, we take a connected graph as a topological space such that the first Betti number of the graph is greater than $1$ and we observe the result in the digital setting. Our next example is the digital interpretation of the topological complexity number of the wedge of $k-$spheres. Finally, we deal with a compact orientable surface of genus $g$. After giving a counterexample in digital images, we reveal how the digital topological complexity numbers of the digital simple closed surface with genus $g$ works in digital images. \\

\acknowledgment{The first author is granted as fellowship by the Scientific and Technological Research Council of Turkey TUBITAK-2211-A. In addition, this work was partially supported by Research Fund of the Ege University (Project Number: FDK-2020-21123)}


\begin{thebibliography}{99}

\bibitem{BorVer:2018} Borat A, Vergili T. Digital lusternik-schnirelmann category. Turkish Journal of Mathematics 2018; 42: 1845-1852.	

\bibitem{Boxer:1999} Boxer L. A classical construction for the digital fundamental group. Journal of Mathematical Imaging and Vision 1999; 10: 51-62.

\bibitem{Boxer:2005} Boxer L. Properties of digital homotopy. Journal of Mathematical Imaging and Vision 2005; 22: 19-26.

\bibitem{Boxer:2006} Boxer L. Homotopy properties of sphere-like digital images. Journal of Mathematical Imaging and Vision 2006; 24: 167-175.

\bibitem{Boxer2:2006} Boxer L. Digital products, wedges, and covering spaces. Journal of Mathematical Imaging and Vision 2006; 25: 169-171.

\bibitem{BoxKar:2012} Boxer L, Karaca, I. Fundemental groups for digital products. Advances and Applications in Mathematical Sciences 2012; 11(4): 161-180.

\bibitem{Gulseli:2014} Burak G. Dijital Kohomoloji Gruplar{\i}. PhD, Adnan Menderes University, Ayd{\i}n, Turkey, 2014.

\bibitem{GulseliKaraca:2017} Burak G, Karaca I. Digital Borsuk-Ulam theorem. Bulletin of the Iranian Mathematical Society 2017; 43: 477-499.

\bibitem{Chen:2004} Chen L. Discrete surfaces and manifolds: A theory of digital-
discrete geometry and topology. Scientific \& Practical Computing, Rockville, MD, 2004.

\bibitem{ChenRong:2010} Chen L, Rong Y. Digital topological method for computing genus and the Betti numbers. Topology and its Applications 2010; 157(12): 1931-1936.

\bibitem{Davis:2017} Davis D. The symmetrized topological complexity of the circle. New York Journal of Mathematics. 2017; 23: 593-602.

\bibitem{Davis2:2017} Davis D. An Approach to the Topological Complexity of the Klein Bottle. (Preprint arXiv: 1612.02747v5)

\bibitem{EgeKaraca:2013} Ege O, Karaca I. Cohomology theory for digital images. Romanian Journal of Information Science and Technology 2013; 16(1): 10-28.

\bibitem{Farber:2003} Farber M. Topological complexity of motion planning. Discrete and Computational Geometry. 2003; 29: 211-221.

\bibitem{Farber:2004} Farber M. Instabilities of Robot Motion. Topology and its Applications 2004; 140: 245-266.

\bibitem{Farber:2006} Farber M. Topology of Robot Motion Planning, In 'Morse Theoretic Methods in Nonlinear Analysis and in Symplectic Topology'. Paul Biran, Octav Cornea, Francois Lalonde editors Springer 2006; pp. 185-230.

\bibitem{Farber:2008} Farber M. Invitation to topological robotics. Zur. Lect. Adv. Math. EMS, 2008.

\bibitem{GrantLuptonOprea:2013} Grant M, Lupton, G, Oprea J. Spaces of topological complexity one. Homology Homotpy and Applications 2013; 15(2): 73-81.

\bibitem{Han:2007} Han SE. Digital fundamental group and Euler characteristic of a connected sum of digital closed surfaces. Information Sciences 2007; 177(16): 3314-3326.

\bibitem{Han:2005} Han SE. Non-Product property of the digital fundamental group. Information Sciences 2005; 171: 73-91.

\bibitem{Herman:1993} Herman GT. Oriented surfaces in digital spaces. CVGIP: Graphical models and image processing 1993; 55: 381-396.

\bibitem{KaracaIs:2018} Karaca I, Is M. Digital topological complexity numbers. Turkish Journal of Mathematics 2018; 42(6): 3173-3181.

\bibitem{MelihKaraca} Is M, Karaca I. The higher topological complexity in digital images (Preprint).

\bibitem{MelihKaraca2} Is M, Karaca I. Certain Topological Methods For Computing Higher Topological Complexity In Digital Images (Preprint).

\bibitem{MorgenRosen:1981} Morgenthaler DG, Rosenfeld, A. Surfaces in three-dimensional images. Informational Control 1981; 51: 227-247.

\bibitem{Ros:1979} Rosenfeld A. Digital topology. American Mathematical Monthly 1979; 86: 76-87.

\bibitem{Ros:1970} Rosenfeld A. Connectivity in digital pictures. Journal of the ACM 1970; 17: 146-160.

\bibitem{Rudyak:2010} Rudyak Y. On higher analogs of topological complexity. Topology and Its Applications 2010; 157(5): 916-920.

\end{thebibliography}
\end{document}